\documentclass{amsart}
\usepackage{amsmath}
\usepackage{amsfonts}
\usepackage{amssymb}
\usepackage{graphicx}
\usepackage{epstopdf}
\usepackage{url}
\usepackage[section]{placeins}
\usepackage[T1]{fontenc}
\usepackage[latin1]{inputenc}

\begin{document}
\title {A note on the possibility of proving the Riemann Hypothesis}

\author{Michele Fanelli}
\email[Michele Fanelli]{michele31fanelli@gmail.com}
\address [Michele Fanelli]{Via L.B. Alberti 5, 20149 Milano, Italy}
\author{Alberto Fanelli}
\email[Alberto Fanelli]{apuntoeffe@gmail.com}
\address[Alberto Fanelli] {Via L.B. Alberti 12, 20149 Milano, Italy}
\date{Feb 15, 2019}

%
%
\begin{abstract}
As well known, the important hypothesis formulated by B.G. RIEMANN in 1859 states that all non-trivial zeroes of the Zeta function $Z(s)=\sum_{n=1}^{\infty } n^{-s}$ should fall on the Critical Line (C.L.) $Re(s)=\frac{1}{2}$.\\
Although direct numerical search of the zeroes failed to identify any outlier, i.e. any zeroes with $Re(s)\neq\frac{1}{2}$, a general proof of the Hypothesis has not yet been found.\\
The present Note aims to approach the problem from a 'reductio ad absurdum' way, i.e. it assumes that an outlier pair of c.c. zero-points, $s=\frac{1}{2}+\xi \pm i.Y_H$ with $\xi\neq0$, has been found, and then proceeds to analyze what are the implications of this assumption. 
Starting from the well-known GRAM-BACKLUND formula for an explicit expression of the Zeta function, the Fundamental Theorem of Algebra (FTA) allows to evidence, through legitimate algebraic manipulations, the necessity that the assumed outlying pair of c.c. zero-points fulfils an implicit additional constraint. The zero-condition according to the GRAM- BACKLUND  formulation and this additional constraint are seen to be mutually incompatible unless  the pair of  c.c. zero-points belong to the C.L.. This conclusion is equivalent to a verification of the RIEMANN Hypothesis.
\end{abstract}

\maketitle

%
%

%
%
\section{Assuming the falsity of the R.H.: the necessary implications}
Since it is known that the original expression of the Zeta Function, $Z(s)=\sum_{n=1}^{\infty } n^{-s}$, leads to a non-convergent series for $Re(s)<1$, and it is known that all non-trivial zeroes of the Zeta Function are contained inside the Critical Strip (C.S.) $0<Re(s)<1$, in order to investigate the conditions to be fulfilled by the zero-points it is necessary to work on a convergent extension of $Z(s)$.\\ For our aims it has proved convenient to choose the GRAM-BACKLUND extension (see \cite{Gram03} and \cite{Back18}) as quoted by EDWARDS (see \cite{Edwa74}), according to which when the zero-condition of the Zeta Function is fulfilled at $s=\frac{1}{2}+\xi \pm i.Y_H$ an equivalent condition, obtained by applying to the Zeta Function the EULER-McLAURIN summation by integration, is necessarily satisfied, to wit: 

\begin{equation}
\text{if } Z(s)=0 \text{ then } Z^*(s)=s.N^{1-s}.[\frac{1}{s.(s-1)}+\frac{1}{Q(s)}]=0
\label{1}
\end{equation}
where
\begin{equation*}
\begin{split}
\frac{1}{Q(s)}=&\{\underset{\mu =2}{\overset{\nu }{\sum }}\frac{B_{2\mu }}{(2\mu )!}.N^{-2\mu }.(s+1).(s+2)\ldots (s+2.\mu -2)\\
&+\frac{1}{s}[\frac{1}{N}.\left(\underset{n=1}{\overset{N-1}{\sum }}\left(\frac{n}{N}\right)^{-s}\right)+\frac{1}{2}+N^s.R_{2\nu }]\}
\end{split}
\end{equation*}
where N is a large integer, the coefficients $B_{2\mu }$ are BERNOULLI numbers and the residual error term $R_{2\nu }$ can be made as small as small as required by a suitable choice of $N$ and $\nu$. 
The considerations developed in the present Section aim to analyze the consistency of the implications of the assumed existence of an outlying c.c. pair of zeroes of the Zeta Function with the specific constraints to which any such pair must answer. \\
To this end, it is necessary to make explicit:

\begin{itemize}
 \item \textbf{The distance of the assumed outlying zero-pair from the C.L..}

This distance ($\xi$ in the following ) is zero if the assumed outlier cannot exist, as is the case of the 'canonical' zero-pairs, and must fall within the bonds $0<|\xi |<\frac{1}{2}$ if $\xi\neq0$ for any hypothetical outlier.
 \item \textbf{The condition occurring when the complex coordinate $s=\frac{1}{2}+\xi \pm i.Y_H$ of a trial point scanning the critical strip falls hypothetically on a zero-point of the Zeta Function.}

Following our choice of the GRAM-BACKLUNG extension of the Zeta Function, in any zero-point $s=\frac{1}{2}+\xi \pm i.Y_H$ in which $Z(s)=0$ the zero-condition is taken to be, see Eq.\eqref{1}: $Z^*(s)=s.N^{1-s}.[\frac{1}{s.(s-1)}+\frac{1}{Q(s)}]=0$, i.e. $\frac{1}{s.(s-1)}+\frac{1}{Q(s)}=0$. The factor $s.N^{1-s}$ in Eq.\eqref{1} is different from zero, and furthermore it is known that $Z(0)\neq0$ and $Z(1)\neq0$. Therefore we can legitimately proceed taking for our working expression of the zero-condition (used to investigate the necessary implications of the zero-condition) the following form:
\begin{equation*}
s.(s-1)+Q(s)=0
\end{equation*}
\end{itemize}

%
%
\section{The implications of the Fundamental Theorem of Algebra (FTA)}
\begin{figure}[!htb]
\centering
\includegraphics[scale=0.49]{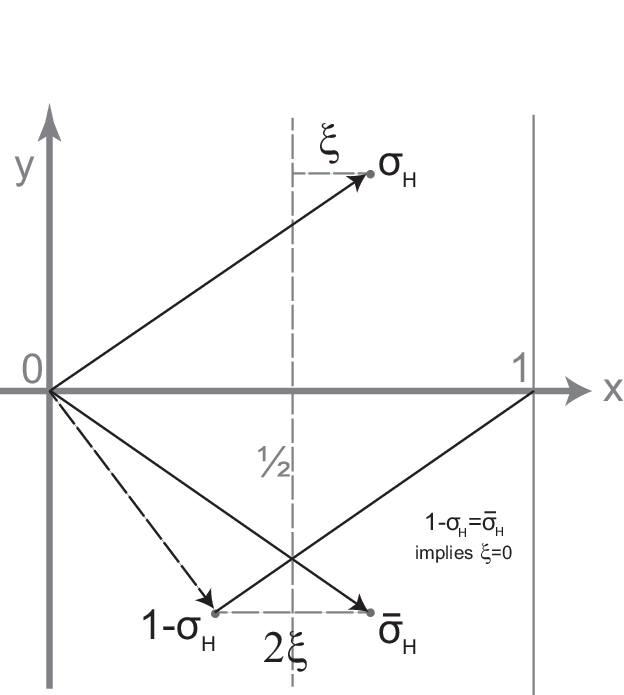}
\caption{Geometric representation of the inequality $1-\sigma_H\neq \overline{\sigma}_H$ for $\xi \neq 0$}
\label{fig:1}
\end{figure}

It is known from FTA that any polynomial equation of the type $P(s)=0$ of degree $N_P>1$ can be factorized as $P(s)=\underset{m=1}{\overset{N_P}{\prod }}\left(s-\sigma _m\right)=0$ where $\sigma _m$ are the zero-points of the polynomial $P(s)$.\\
For our purpose the assumed c.c. pair of zero-points is henceforth denoted as $\sigma _H=\frac{1}{2}+\xi +i.Y_H$ and $\overline{\sigma }_H=\frac{1}{2}+\xi -i.Y_H$.
The zero condition for the Zeta function implies necessarily $Z^*(s)=0$ (see Eq.\eqref{1}), which can be treated as a polynomial identity when $s=\sigma _H$ or $s=\overline{\sigma}_H$ because $Q(\sigma_H)$ is implicitly assumed to be a known complex quantity.\\
This implication of the zero-condition, namely $s.(s-1)+Q(s)=s^2-s+Q(s)=0$, is assumedly fulfilled by $s=\sigma _H=\frac{1}{2}+\xi +i.Y_H$ and by $s=\overline{\sigma}_H=\frac{1}{2}+\xi -i.Y_H$.\\
Therefore it comes:
\begin{equation}
\begin{split}
\sigma _H.(\sigma_H-1)+Q(\sigma _H)&\equiv 0 \\
\text{ and } \overline{\sigma}_H.(\overline{\sigma}_H-1)+Q(\overline{\sigma}_H)&\equiv 0
\end{split}
\label{2}
\end{equation}
After FTA, the implication of the pseudo-polynomial zero-condition, i.e. $s.(s-1)+Q(\sigma_H )=0$, can be factorized as: 
\begin{equation} 
(s-\sigma_H).(s-\overline{\sigma}_H)=0
\label{3}
\end{equation}
On the other hand, from the polynomial division $\frac{s.(s-1)+Q(\sigma_H)}{s-\sigma_H}$ one gets:
\begin{equation*}
\frac{s.(s-1)+Q(\sigma_H)}{s-\sigma_H}=s-(1-\sigma_H)
\end{equation*}
and likewise from the pseudo-polynomial division $\frac{s.(s-1)+Q(\overline{\sigma}_H)}{s-\overline{\sigma }_H}$ one gets:
\begin{equation*}
\frac{s.(s-1)+Q(\overline{\sigma}_H)}{s-\overline{\sigma}_H}=s-(1-\overline{\sigma}_H)
\end{equation*}
with rests $s.(\sigma _H-1)+Q(\sigma _H)$ or, respectively, $s.(\overline{\sigma }_H-1)+Q(\overline{\sigma }_H)$, which according to our assumptions are zero for $s=\sigma_H$ or $s=\overline{\sigma}_H$, see Eq.\eqref{2}, so that it comes: 
\begin{equation*}
s.(s-1)+Q(\sigma _H)=(s-\sigma _H).[s-(1-\sigma _H)] \text{ and } s.(s-1)+Q(\overline{\sigma }_H)=(s-\overline{\sigma }_H).[s-(1-\overline{\sigma }_H)]
\end{equation*}
Therefore our assumptions lead to two equally necessary implications of the zero-condition, namely:
\begin{equation*}
(s-\sigma _H).(s-\overline{\sigma }_H)=0 \text{ and } (s-\sigma _H).[s-(1-\sigma _H)]=0
\end{equation*}
from which, for $s=\overline{\sigma}_H$  :
\begin{equation*}
\overline{\sigma }_H=1-\sigma _H  \text{ (see note} \footnote{As a check, let us work out symbolically  the product  $(s-\sigma_H ).[s-(1-\sigma_H )]$, keeping account that   according to our assumptions   $\sigma_H^2-\sigma_H+Q(\sigma_H )=0$,  or  $-(\sigma_H^2-\sigma_H )=Q(\sigma_H )$; it follows that
 $(s-\sigma_H ).[s-(1-\sigma_H )]=(s-\sigma_H ).(s+\sigma_H-1)=s^2-\sigma_H^2-s+\sigma_H=s^2-s-(\sigma_H^2-\sigma_H )=s^2-s+Q(\sigma_H )=s^2-s+Q(\sigma_H )$,     so that if   $s\neq 1-\sigma_H$      it must be     $s=\sigma_H$;   likewise    if   $s\neq\sigma_H$     it must be     $s=1-\sigma_H$. On the other hand if  $s\neq\sigma_H$,   according to our assumptions   it must be    $s=\overline{\sigma}_H$,    and  the two conditions  are compatible only if  $1-\sigma_H=\overline{\sigma}_H$.   Likewise  we find that   if        $s\neq\sigma_H$      it must be    $s=1-\sigma_H$,   but   according to our assumptions   it must be    $s=\overline{\sigma}_H$, and  again the two conditions are compatible only if  $\overline{\sigma}_H=1-\sigma_H$.}  \text{)}
\end{equation*}
which implies, after the assumed definition $\sigma _H=\frac{1}{2}+\xi +i.Y_H$:
\begin{equation*}
\overline{\sigma }_H=\frac{1}{2}+\xi -i.Y_H=1-\sigma _H=\frac{1}{2}-\xi -i.Y_H
\end{equation*}
which cannot be fulfilled unless
\begin{equation*}
{\xi=0}
\end{equation*}
thus contradicting our assumption $\xi\neq 0$ (see Figure ~\ref{fig:1}). 

%
%
\section{Corollary}
From Eq.\eqref{3}, i.e. $(s-\sigma_H).(s-\overline{\sigma}_H)=0$, or $s^2-(\sigma _H+\overline{\sigma }_H).s+\sigma _H.\overline{\sigma }_H=0$, which must be equivalent to $s^2-s+Q(\sigma_H)=0$, it comes:
\begin{equation*}
-(\sigma_H+\overline{\sigma}_H)=-(1+2.\xi )=-1
\end{equation*}
 which confirms that it must be $\xi=0$ and $Q(\sigma_H)=\sigma_H.\overline{\sigma}_H=\frac{1}{4}+\xi ^2+\xi +Y_H^2$ from which:
$Q\left(\sigma _H\right)\in \mathbb{R}$, and since $\xi=0$ it comes $Q(\sigma_H)=Q(\overline{\sigma}_H)=\frac{1}{4}+Y_H^2$

%
%
\section{Conclusions}
The considerations developed in the preceding section of the present note seem to show that the zero-condition for the Zeta Function based on the GRAM-BACKLUND extension of the Zeta Function would necessarily entail an internal inconsistency unless $Re(s)=\frac{1}{2}$, i.e. unless the R.H. is verified.\\
This appears to depend on the invariance of the product $s.(1-s)$ upon the exchange $s\leftrightarrow(1-s)$.\\
It looks difficult to put in doubt the above conclusion, unless some fundamental flaw of the GRAM-BACKLUND representation of the Zeta Function could be identified. The extensive investigations about the actual zero-points of the Zeta Function carried out in the past, see ODLYZKO \cite{Odly01}, would seem, however, to rule out such an eventuality.\\
On the other hand, if the present Note should be found without flaw, it remains to understand why the approach outlined in the present Essay was left until now unexplored and neglected.

%
\section{Figures}

\end{document}